\documentclass[a4paper,10pt]{amsart}
\usepackage{amsmath,latexsym,amssymb,amsthm}
\usepackage[dvips]{graphicx}
\usepackage[english]{babel}

\title{Counting and computing regions of $D$-decomposition: algebro-geometric approach.}
\author{Oleg O. Vasil'ev}
\address{Ya.Z. Tsypkin Laboratory ``Adaptive and Robust Control Systems'', V.A.Trapesnikov Institute of Control Sciences of Russian Academy of Sciences, Moscow, Russian Federation;\newline Chair of Applied Mathematics and Computer Modelling, Department of Automation and Computer Science,  I.M. Gubkin Russian State University of Oil and Gas, Moscow, Russian Federation.}
\email{ovasilev@gmail.com}
\thanks{The author was partially supported by RFBR grant 11-08-00223-à}
\begin{document}

\newcommand {\R}{\mathbb R}
\newcommand {\p}{\mathfrak p}
\newcommand {\Z}{\mathbb Z}
\newcommand {\N}{\mathbb N}
\newcommand {\CN}{\mathbb C}
\newcommand {\K}{\mathbb K}
\newcommand {\F}{\mathbb F}
\newcommand {\B}{\mathbb B}
\newcommand {\NN}{\mathbb N\cup\{0\}}
\theoremstyle{definition}
\newtheorem{definition}{Definition}
\theoremstyle{plain}
\newtheorem{theorem}{Theorem}
\newtheorem{theor}{Òåîðåìà}
\newtheorem{lemma}{Lemma}
\newtheorem{corollary}{Corollary}
\newtheorem{rem}{Note}
\newtheorem{example}{Example}
\begin{abstract}
New methods for $D$-decomposition analysis are presented. They are based on topology of real algebraic varieties and computational real algebraic geometry. The estimate of  number of root invariant regions for polynomial parametric families of polynomial and matrices is given. For the case of two parametric family more sharp estimate is proven. Theoretic results are supported by various numerical simulations that show higher precision of presented methods with respect to traditional ones.
The presented methods are inherently global and could be applied for studying $D$-decomposition for the space of parameters as a whole instead of some prescribed regions.
For symbolic computations the Maple v.14 software and its package RegularChains are used.
\end{abstract}
\maketitle

\section{Introduction}
In this paper linear differential equation with constant coefficients 
$$
\dot{x}(t)=Ax(t),
$$
where $A$ is a complex(real) $n\times n$ matrix, will be called linear continuous system.

Similarly, linear discrete system is a finite-difference equation
$$
x(t+1)=Ax(t),
$$
where $A$ is a complex(real) $n\times n$ matrix.

The problem of finding stability domains of a linear control system depending on a vector parameter is a classical one. 
It can be formulated as such:
Let us consider a polynomial  $P(s,k)$ with  complex(for example, real) coefficients. polynomially depending on real parameters $k=(k_1,\ldots,k_l)\in \K.$  Complex parameters  could be considered as  $a+jb,$ where $a, b$ are real. 

This polynomial is usually considered as a characteristic polynomial of a linear system. Vector $k$ understood as indeterminacy or parameters of a regulator.

We need to find such regions $U\subseteq \K$ that for all $k\in U$ $P(s,k)$ will be continuous-time(resp. discrete-time) stable.

Polynomial is called continuous-time stable(Hurwitz) if all of its roots have negative real part. Polynomial is called discrete-time stable(Schur) if all of it's roots is have modulus lesser than $1$.

Continuous-time(res. discrete-time) system is stable if and only if its characteristic polynomial is a Hurwitz(resp. Schur) one.

It is obvious that the problems of finding regions of hurwitzness or schurness of a polynomial is equivalent. Clearly, polynomial $P(z)$ is Schur if and only if $(s-1)^{\deg P}P(\frac{s+1}{s-1})$ is Hurwitz. Therefore we can consider only continuous-time systems, and we will do that in any case, unless otherwise specified.

There are two cases which are of the greater interest for control theoretist, namely the case of characteristic polynomial of SISO system with affine dependency on parameters and the case of family of matrices affinely depending on parameters,  for example, in the case of construction of a linear output feedback.  In the latter case polynomial $P(s,k),$ will be $\det (sI -(A_0+\sum_{i=1}^lk_iA_i)).$

One of a methods of studying stability regions is a $D$-decomposition approach based on study of decomposition of parameter space into connected regions with equal numbers of stable and unstable roots(root-invariant regions or regions of $D$-decomposition). It is clear that the border of these regions is a semialgebraic hypersurface in $K$ parametrized by equalities $\Re(P(j\omega,k))=0, \Im (P(j\omega,k))=0$ and, possibly, a hypersurface with an equation  $|a_{\deg P}(k)|=0,$ where $a_{\deg P}$ \-- is a leading term of $P.$ $\Re$ and $\Im$ here denote real and imaginary part of an expression.

The main idea of this approach goes back to I. Vishnegradsky\cite{Vish1876}, but first clear formulation and detailed study of a method have been done by Yu.I. Neimark \cite{Nei1948}, \cite{Nei1949}. 
The contemporary state-of-the-art could be obtained from paper  \cite{GP2006} and review \cite{GPT2008}.

The main idea of our paper is an application of various methods of computational algebra and real algebraic geometry to the study of $D$-decomposition\-- both computational(finding sample points from all root invariant regions) and theoretical \-- upper bounds to the number of regions. 

The famous Tarski-Seidenberg theorem (\cite{Tar1951,Sei1954}) tells us that the problem of description of all root-invariant regions, as well as stable ones is algorithmically solvable.  Nevertheless, effective algorithms of study of real polynomial systems of equations have been  included in general purpose computer algebra systems only recently. For example, Maple package $RegularChains$ was included in Maple only in 2005 and it is still rapidly developing. The history of those algorithms itself goes back to discovery of cylindric algebraic decomposition by G.E. Collins in 1975 \cite{Col1975} and partial cylindric algebraic decomposition by G.E.Collins and H.Hong in 1991\cite{CH1991}. Monograph \cite{BPC2006} is a standart reference book on the subject.

Moreover, applications of quantifier elimination algorithms to control theory goes back to the seminal paper by B.D.O. Anderson, N.K. Bose, E.I.Jury published in 1975 \cite{ABJ1975}, where an application of quantifier elimination to the linear output feedback stabilisation problem have been studied.  They used Routh-Hurwitz theory, multivariate resultants and classic Tarski-Seidenberg quantifier elimination technique to check existence or non-existence of stable output feedback. However, there was no effective computational techniques, software and sufficiently powerful machines at the time, so the interest to this topic have risen only in the midde of 1990's\cite{DYA1997, Neu1997, Jir1996}. The application of real algebraic geometry to the study of the stable points set and the Nyquist map is stated as an unsolved problem in monograph by E.A. Jonckheere \cite{Jon1997} The approach used in this papers is quite similar to the Anderson-one, but authors have used not an original Tarski-Seidenberg algorithms but  cylindrical algebraic decomposition implemented at that time in the QEPCAD software by H.Hong.  His collaborative work with R. Liska and S.Steinberg \cite{HLS1996} gave rise to applications of partial cylindric algebraic decomposition algorithms to stability problems for PDE. 

Later, H.Anai and S.Hara et al. in series of papers have introduced the new methods of robust control synthesis with quantifier elimination based on so-called Sturm-Habicht sequence, Symbolic-Numeric CAD etc (see \cite{AH2006},\cite{IYAY2009} and references therein), they have done computations with these methods in computer algebra system Risa/Asir. They have also  designed a Maple package SyNRAC, where some of their algorithms was implemented.

Some similar considerations in the context of Nyquist criterium could be found in papers by  N.P. Ke (\cite{Ke2000} and references therein). 
 The paper \cite{SXZ2011} is devoted to the applications of computational real algebraic geometry to an asymptotic stability analysis of nonlinear systems. A.D. Bruno et al. \cite{BBV2012} introduces use of resultants for study of $D$-decomposition border.

This paper is organised as follows. In section 2 we obtain some upper bounds to the number of root invariant regions in the parameter space in two-dimensional and general cases. Using both Groebner bases and classical resultant-based elimination technique we obtain a variety containing projection of intersection of hypersurfaces $\Re(P(j\omega,k))=0, \Im (P(j\omega,k))=0$ onto $\K.$
Then we estimate the number of regions of its complement using classical Harnack and Bezout theorems  in two-dimensional case and Warren inequality in the general case. These results are far-going generalisation of upper bounds from \cite{Gry2004, GP2006}.

In section 3 we give the description of Maple-based algorithm of studying $D$-decomposition and concentrate on some particular examples from \cite{GP2006}, which illustrate main features of our approach. We use a combination of Groebner basis-based or resultant-based elimination and partial cylindric algebraic decomposition. Specifically, for the continuous-time system first appeared in \cite{QBRDYD1995}, Example 2 and studied in \cite{GP2006}, Example 15, we have found one new root-invariant region. This example tell us about one very important feature of our methodology \-- its inherent globality. We automatically obtain the structure of $D$-decomposition not only for some predefined small region, but for parameter space as a whole.
\section{Estimates for the number of root invariant regions in the parameter space}

Let us consider a polynomial  $P(s,k)$ with  complex(for example, real) coefficients. polynomially depending on real parameters $k=(k_1,\ldots,k_l)\in \K.$  Complex parameters  could be considered as  $a+jb,$ where $a, b$ are real. 

This polynomial is usually considered as a characteristic polynomial of a linear system. Vector $k$ understood as indeterminacy or parameters of a regulator.

We need to find such regions $U\subseteq \K$ that for all $k\in U$ $P(s,k)$ will be continuous-time(resp. discrete-time) stable. The case of discrete-time stability is equivalent to continous-time one up to linear fractional transformation. Therefore we can consider only the continous-time case unlessotherwise is specified.

The $D$-decomposition approach is based on study of decomposition of parameter space into connected regions with equal numbers of stable and unstable roots(root-invariant regions or regions of $D$-decomposition). 

It is clear that the border of these regions is a semialgebraic hypersurface in $K$ parametrized by equalities $\Re(P(j\omega,k))=0, \Im (P(j\omega,k))=0$ and, possibly, a hypersurface with an equation  $|a_{\deg P}(k)|=0,$ where $a_{\deg P}$ \-- is a leading term of $P.$ $\Re$ and $\Im$ here denote real and imaginary part of an expression.

In the following theorems there are obtained  upper bounds for a number of root invariant regions.

\begin{theorem}\label{th:1}
Let $P(s,k_1,\ldots,k_n)$  be a complex(in particular, may be real) polynomial of degree $t$ on $s$ and of degree $d$ on all $k_i$ together.
Suppose that 
$\Re( P(j\omega,k_1,\ldots,k_n))$ and $\Im( P(j\omega,k_1,\ldots,k_n)),$ do not have any common divisors non-trivially depending on $\omega.$

Then the number of regions of $D$-decomposition is no greater than
$$
6(4td+4d)^n,
$$
\end{theorem}
\begin{proof}

The border of $D$-decomposition in the continuous-time case is a semialgebraic set in the parameter space $p_1,\ldots,p_n$.
It is defined parametrically by the pair of equations  $\Re(P(j\omega,k_1,\ldots,k_n))=0$, $\Im(P(j\omega,k_1,\ldots,k_n))=0,$ and an equation $a_t(k_1,\ldots,k_n)=0.$ $a_t$ is a coefficient of  $s^t$ in $P(s,k_1,\ldots,k_n)$(\cite{PS2002} \S 4.1.2).

This semialgebraic set could be continued to a minimal algebraic variety containing it. Because $\Re(P(j\omega,k_1,\ldots,k_n), \Im(P(j\omega,k_1,\ldots,k_n))$ are coprime over $\CN[\omega]$  its dimension will stay the same, therefore the number of connected components of complement will not be lesser.

Consider a minimal algebraic variety $X$ containing the set parametrized by $\Re(P(j\omega,k_1,\ldots,k_n))=0, \Im(P(j\omega,k_1,\ldots,k_n))=0.$ Let $g$ be a Groebner basis of an ideal generated by parametrization polynomials relative to some $\omega$-eliminating order. Then by an elimination theorem (\cite{KOL2000}, Ch.3 \S 1, Theorem 2; \cite{KOL2000}, Ch.3 \S 1, Ex.5) $X$ is defined by elements of $g$ that do not contain $\omega.$

By Proposition 1 in  \cite{KOL2000}, Ch.3 \S 6 the resultant $r$ of $\Re(P(j\omega,k_1,\ldots,k_n))$, $\Im(P(j\omega,k_1,\ldots,k_n))$ in variable $\omega$ defines an algebraic set containing the border of $D$-decomposition.

This set could be bigger than the set defined by a basis. But because of coprimeness of $\Re(P(j\omega,k_1,\ldots,k_n)), \Im(P(j\omega,k_1,\ldots,k_n))$ up to $\omega$ resultant $r$ is a non-zero polynomial. Therefore the set of zeros of $r$ divides $\K$ in a number of parts that is no greater than the number of parts generated by zeros of the Groebner basis. 

In order to find that number we need to estimate the degree of $h=ra_t$(or $h=r(\Re^2(a_t)+\Im^2(a_t))$)  It can be easily shown that it is no greater that $2td+2d.$
Finally, we need to bound the number of regions. Using Theorem 2 \cite{War1968} we get $6(2 \deg h)^n.$

The bound for continuous-time is the same as for discrete-time one because linear fractional transformation does not change the degree of polynomials.
\end{proof}

It must be noted that the main result of \cite{BPC2009} could signify that the multiplier $6$ in formulation of the theorem above could be replaced by $2.$

In $2$-dimensional case we can obtain better bound, even sharp in some sense.
To prove an upper bound for the number of root invariant regions we have to prove the following lemma:
\begin{lemma}\label{lem:1}
Let $X$ be the plane real (possibly singular and reducible) affine algebraic curve of degree $n.$

Then its complement in $\R^2$ consists of no more than  $\frac{n^2+n+2}{2}$ connected components.

This bound is sharp and it is reached on union of lines in general position.
\end{lemma}
\begin{proof}
If $X$ is a real irreducible nonsingular projective curve of degree $n$ then, by Harnack inequality, (\cite{Gud1974} \S 3, \cite{JLC1931} Ch.4 Th. 11) it divides plane to no more than $\frac{(n-1)(n-2)}{2}+2$ components and $\frac{(n-1)(n-2)}{2}+1$ of them correspond to the ovals of $X.$

If $X$ is singular, then by First Harnack Theorem  (\cite{JLC1931}, Ch.4, Th.12) real irreducible projective curve has no more than $\frac{(n-1)(n-2)}{2}-\sum k_i(k_i-1)+1$ branches. Here $k_i$  are the multiplicities of singular points.

Note that every self-intersection or self-tangency of $l$ fragments of ovals of the curve in a singular point gives birth to no more than $2l$ new connected regions in a small neiborghood of a point. All other types of singularities do not lead to an increase of number of regions.

Therefore if the mutiplicity of a singular point has increased by $l,$ then the maximal possible number of ovals is reduced by $(k_i+l)(k_i+l-1).$ Hence the number of connected components of a complement to $X$ increases by no more than $2l-(k_i+l)(k_i+l-1),$ but if  $k_i, l \geq 1, $ then $2l-(k_i+l)(k_i+l-1)\leq 0.$  Thus singular irreducible projective curve divides the projective plane into no more than $\frac{(n-1)(n-2)}{2}+2$ components.

Let $X$ be a reducible curve. Denote by the symbols  $m\vdash n$ such a fact that $m=(m_1,\ldots,m_{|m|})$ is a partition of $n$ to the degrees of irreducible components of $X.$

Let  $k=|\{i|m_i=1\}|$, $l=|\{i|m_i>1\}|.$

We can estimate the number of regions $Reg$ as
$$
\begin{array}{l}
Reg\leq \max_{m\vdash n}(\sum_{m_i\in m}\frac{(m_i-1)(m_i-2)+2}{2}+\sum_{i\neq j}m_im_j+n-k+1)=\\=\frac{n^2-n+2}{2}+\max_{m\vdash n}(2|m|-k)=\frac{n^2-n+2}{2}+\max_{m\vdash n}(2l+k)\leq\\\leq\frac{n^2-n+2}{2}+\max_{2l+k\leq n, 0\leq l,0\leq k}(2l+k)=\frac{n^2+n+2}{2}.
\end{array}
$$

This estimate is obviously sharp. It is sufficient to take a reducible curve that consists from $n$ lines in general position(e.g. without parallels or more than double intersections).

\begin{figure}\label{locstr}
\includegraphics[scale=0.19]{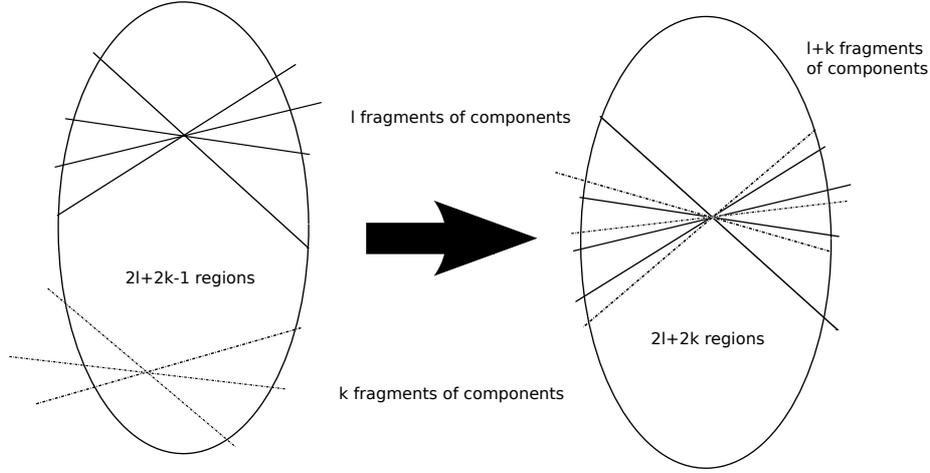}
\caption{Local structure self-intersection of a reducible and/or singular curve in a neiboghrood of self-interection point up to isotopy}
\end{figure}

The only thing we have to prove is the following inequality:
$$
Reg\leq \max_{m\vdash n}(\sum_{m_i\in m}\frac{(m_i-1)(m_i-2)+2}{2}+\sum_{i\neq j}m_im_j+n-k+1)
$$
The first term in the left part of inequality is trivially derived from previous considerations. The term $n-k$ corresponds to the  transition from the projective case to affine one. Actually, there are no more than $n$ intersections between $X$ and the infinite line. $k$ of them are intersections between linear components and the infinite line, but this intersections do not add any new connected components.

We can deduce from Bezout theorem that there are no more than $m_im_j$ intersections between irreducible components of degrees $m_i$, $m_j.$ Therefore we need to prove that every such intersection creates no more connected components than its multiplicity.

Let us consider sufficiently small neighborhoods  of points of two intersecting and non-intersecting irreducible components(see Figure 1).  Those irreducible components are locally isotopic to an intersection of lines(one line if non-singular), therefore we need to consider only the case of lines , which is shown at the figure 1.

Now we need globalize the local structure. Let us consider the following numeration of regions.

Suppose that we have $p$ irreducible components. Let us denote by $N(i)$ a number of regions that the irreducible component $i$ adds. 

Let us divide the naturals into $p+1$ parts. \-- from $2$ to $2+N(1)-1$, from $2+N(1)$ to $2+N(1)+N(2)-1$, from $N=2+\sum_{i=1}^pN(i)$ to $+\infty.$
We can break the plane into the finite union of neighborhoods such that each of them contains no more than one point of intersection of components. 

We can enumerate connected regions appearing in those neighborhoods.

Let us take an arbitrary region that is not internal for any oval of an irreducible component. It will be the first region.

Then we take an arbitrary neighborhood.  If there are a point of intersection of irreducible components in it then we bypass that point counterclockwise and enumerate all appearing non-numbered regions using the following rule:

\begin{enumerate}
 \item If the component $i$ that borders current region from counterclockwise-side have appeared before and all numbers between $2+\sum_{w=0}^{i-1}N(w)$ and $2+\sum_{w=0}^{i}N(w) -1$ already correspond to some region, then let us numerate it by the first free number, that is greater than $N.$

\item If, otherwise, there are free number from $2+\sum_{w=0}^{i-1}N(w)$ to $2+\sum_{w=0}^{i}N(w) -1$ then let us numerate our region by first free number from this range.
\end{enumerate}

Then we use the same rule to enumerate all non-enumerated regions around points of self-intersections of irreducible components. 

Finally, all remaining regions will be enumerated upside-down and if this will not be possible \-- from left to right.

We have $e$ enumerated regions of the plane. 

Let us call  a (possibly non-compact) region on the plane with border containing more than one irreducible component a curvilinear polygon. Points(possibly infinite or even formal\-- i.e. directions in which border component changes) of intersections between two distinct irreducible components will be called vertices of a polygon. Parts of irreducible components lying betwóen vertices will be called edges.

Let $R_1\ldots R_z$ be a curvilinear polygon. Every irreducible component form no more than half of all edges of $R_1\ldots R_z.$

We have established a correspondence between the set of regions and segment of naturals $[1,\ldots, e+1).$
Consider it as a part of real line. Break every interval of form $[i,i+1)$ into parts corresponding to an irreducible components that appear in its border. The length of every subinterval in this partition is proportional to the number of edges formed by  every component. Therefore there will be no element of the partition of length greater than $\frac{1}{2}-\epsilon$ and lesser than $1.$

 Consider a neighborhood of any point of intersection between components. Let us bypass that point in any direction. During the bypass every component would appear to be a border for no more than $2j$ regions, where $j$ is the intersection multiplicity of that component in the point. From the other side the contribution of this component to the partition of interval $[1,\ldots, e+1)$ in the point is no greater than $2j\frac{1}{2}=j.$ This proves the lemma.
\end{proof}
\begin{theorem}\label{th:2}
Let $P(s,r,p)$ be a polynomial of degree $t$ with complex(in particular, possibly, real) coefficients.polynomially depending from two real parameters $r$, $p$ or one complex parameter of a form $r+jp.$ Denote by $d$ its maximal degree as polynomial on variables $r$ and $p$ together.

Suppose that $\Re( P(j\omega,r,p))$ and $\Im( P(j\omega,r,p)),$ as complex polynomials does not have any common divisors  that nontrivially depends on $\omega.$

Then the number of regions of $D$-decomposition is no greater than 
$$
\frac{q^2+q+2}{2},
$$
with $q=2td+2d.$
\end{theorem}
\begin{proof}
Follows from Lemma \ref{lem:1} and the proof of previous theorem.
\end{proof}
This bound could not be generalised onto $n$-dimensional  case, because of lack of any sharp Harnack-like bounds for higher-dimensional real algebraic varieties.
Moreover, we can easily show that the case of hyperplanes in general position even in $3$-dimensional case do not give a maximal number of  regions. 

Namely, Bihan\cite{Bih2003} has shown that the maximal number of connected components of real projective algebraic surface of degree $q$ is equivalent as $q\to\infty$ to $dq^3,$ where $d\in [\frac{13}{36},\frac{5}{12}].$
Moreover  \cite{Vir1998}, from Comessatti-Petrovsky-Oleinik and Smith-Thom inequalities could be obtained an upper bound  $\frac{5}{12}q^3-\frac{3}{2}q^2+\frac{25}{12}q$.

But $q$ planes could divide the space only in $\frac{q^3+5q+6}{6}$ parts.

Now we can write bounds for the case of polynomial matrix family.
\begin{corollary}
Let $A(r,p)$ be a $t\times t$-matrix, which entries are complex(possibly real) polynomials on real parameters $r$,$p$(may be on one complex parameter $r+jp$) of degree no greater than $d$

Suppose that $\Re (\chi_A(j\omega,r,p))$ and $\Im (\chi_A(j\omega,r,p)),$ as complex polynomials does not have any common divisors  that nontrivially depends on $\omega.$

Then the number of regions of $D$-decomposition is no greater than:
$$
\frac{q^2+q+2}{2},
$$
where  $q=2t^2d$ in continuous case and $q=2t^2d+2td$ in discrete one.
\end{corollary}
\begin{corollary}
Let $A(k_1,\ldots,k_n)$ be a $t\times t$-matrix, which entries are complex(possibly real) polynomials on real parameters $k_1\ldots k_n$ of degree no greater than $d$.

Suppose that $\Re( \chi_A(j\omega,k_1,\ldots,k_n))$ and $\Im (\chi_A(j\omega,k_1,\ldots,k_n)),$ as complex polynomials does not have any common divisors  that nontrivially depends on $\omega.$

Then the number of regions of $D$-decomposition in continuous case  is no greater than:
$6(4t^2d)^{n}, $ and no greater than $6(4t^2d+4td)^n$ in discrete one.
\end{corollary}
\section{$D$-decompositions for 2-parametric families: computational study}
In order to study 2-dimensional $D$-decompositions with Maple 14 we have used the following scheme.
\begin{enumerate}
 \item Convertion of family of matrices to the family of polynomials  $P(s,r,p)$, if needed.
 \item If we think about our family as a discrete-time object then we convert it to continuous-time using transformation $P(s)\mapsto (s-1)^{\deg P}P(\frac{s+1}{s-1}).$
 \item Computing an equation defining most of irreducible components of an algebraic curve, containing border of $D$-decomposition using command $EliminationIdeal$ for an ideal generated by $\Re( P(j\omega,r,p))$ and $\Im( P(j\omega,r,p))$ or eliminating $\omega$ using resultants.
\item if $P(s,r,p)$ has a real leading term then we add an irreducible component generated by its leading coefficient $a_{\deg P}(r,p)$ else we add $\Re^2(a_{\deg P}(r,p))+\Im^2(a_{\deg P}(r,p)).$ We get an algebraic curve containing the border of $D$-decomposition.
\item Using  $PartialCylindricalAlgebraicDecomposition,$ we obtain a point cloud that contains at least one point from every region of $D$-decomposition.
\item We count number of stable roots for every point from a cloud and generate the list of stable points.
\end{enumerate}

Let us proceed to examples. 
\begin{example}
Consider a discrete-time object $s^6+(r+jp)s^5+\frac{3}{2}.$ 
Converting it to continuous time we obtain: $p(s,r,p)=(s+1)^6+(r+jp)(s+1)^5(s-1)+\frac{3}{2}(s-1)^6.$

$D$-decomposition of family $z^n+(r+jp)z^{n-1}+\alpha,$ with $\alpha>1$ has $(n-1)^2+1$ regions \cite{GP2006}. In  \cite{GP2006}, there is also given an estimate for number of root invariant regions for families of a form $a(s)+(r+jp)b(s),$ namely $(n-1)^2+2.$ 
Polynomial parametrization of border of $D$-decomposition is:
$$
\begin{array}{l}
 -r-4ps^5+4ps-\frac{75}{2}s^2+5rs^4+5rs^2-rs^6+\frac{5}{2}+\frac{75}{2}s^4-\frac{5}{2}s^6=0\\
-4rs-ps^6-3s^5+5ps^2-p-3s+5ps^4+4rs^5+10s^3=0
\end{array}
$$

Explicit equation of border of $D$-decomposition has degree $10:$
$$
\begin{array}{l}
9216p^{10}+46080p^8r^2+92160p^6r^4+92160p^4r^6+46080p^2r^8+9216r^{10}-\\-94464p^8-377856p^6r^2-566784p^4r^4-377856p^2r^6-94464r^8+301440p^6+\\+683136p^4r^2+1051776p^2r^4+276864r^6-309600p^4-619200p^2r^2-309600r^4+\\+122500p^2+122500r^2-15625=0
\end{array}
$$
By the Lemma \ref{lem:1}, if  degree of a border is $10,$  $D$-decomposition contains no more than $56$ regions. There are  $26$ regions in the example considered. 

There are $56$ regions if and only if the border of $D$-decomposition is an arrangement of $10$ lines without parallels and triple intersections.

At the figure 2 the $D$-decomposition and the point cloud are shown. There are no stability domains.
\begin{figure}\label{4.6}
\includegraphics[scale=0.4]{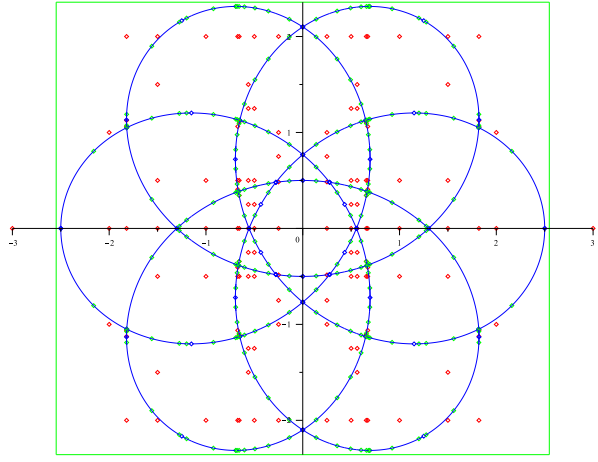}
\caption{$D$-decomposition and a point cloud from connected root invariant regions $s^6+(r+jp)s^5+\frac{3}{2}.$}
\end{figure}
\end{example}
\begin{example}
The other example \cite{GP2006} is a discrete-time object  $z^6+(r+jp)z^5+\frac{3}{20}.$

It is an example of one possible fenomenon appearing in a given way of studying $D$-decomposition.
Namely, an equation  $\Re^2(a_6(r,p))+\Im^2(a_6(r,p))=(r+\frac{23}{20})^2=0$ has nonzero number of roots but it does not separate root invariant regions -- there are the same number of stable ad unstable roots from every side of the line defined by the equation.

The border of $D$-decomposition(without extraneous components) is given by an euation of degree $10:$
$$\begin{array}{l}
9,216\cdot10^{13}p^{10}+4,608\cdot10^{14}p^8r^2+9,216\cdot10^{14}p^6r^4+9,216\cdot10^{15}p^4r^6+\\+4,608\cdot10^{14}p^2r^8+9,216\cdot10^{13}r^{10}+8,1792\cdot10^{13}p^8+3,27168\cdot10^{14}p^6r^2+\\+4,90752\cdot10^{14}p^4r^4+3,27168\cdot10^{14}p^2r^6+8,1792\cdot10^{13}r^8+\\+1,30276416\cdot10^{15}p^6-1,821010752\cdot10^{16}p^4r^2+1,865389248\cdot10^{15}p^2r^4-\\-1,15483584\cdot10^{15}r^6+6,82460784\cdot10^{11}p^4+1,364921568\cdot10^{14}p^2r^2+\\+6,82460784\cdot10^{13}r^4+4,160322828658\cdot10^{15}p^2+4,160322828658\cdot10^{15}r^2-\\-3,573226485213841\cdot10^{15}=0.\end{array}
$$

At the figure 3 the border of $D$-decomposition(including extraneous line)  is shown.
\begin{figure}\label{4.7}
\includegraphics[scale=0.4]{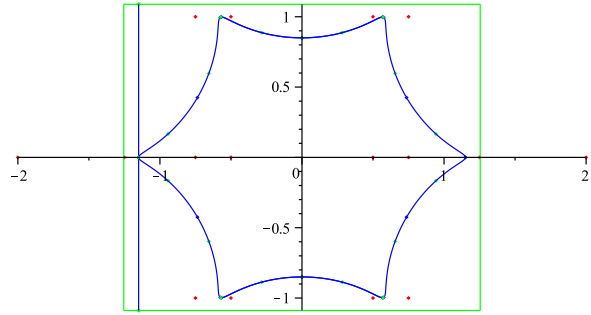}
\caption{The border of $D$-decompotion with extraneous line that corresponds to the case of vanishing of the leading term and  point cloud for the family $z^6+(r+jp)z^5+\frac{3}{20}.$}
\end{figure}
\end{example}
\begin{example}
The other example shows us how a minimal algebraic curve containing the border of $D$-decomposition may differ from the border of $D$-decomposition itself.

Let us consider an output feedback problem of the form $K=\begin{pmatrix}r&p\\-p&r\end{pmatrix}$ for discrete-time system defined by matrices:
$$
A=\begin{pmatrix}0,4753&0,7579&7,9939\\-0,0415&0,8905&0,7579\\-0,0758&-0,0415&0,4753\end{pmatrix}
$$

$$
B=\begin{pmatrix}0,0801&0,0430\\-0,0015&0,0948\\-0,043&-0,0015\end{pmatrix},\qquad C=\begin{pmatrix}1&0&0\\0&1&0\end{pmatrix}
$$

A minimal algebraic curve containing the border of $D$-decomposition consists from three intersecting ovals.
\begin{figure}\label{14}
\includegraphics[scale=0.6]{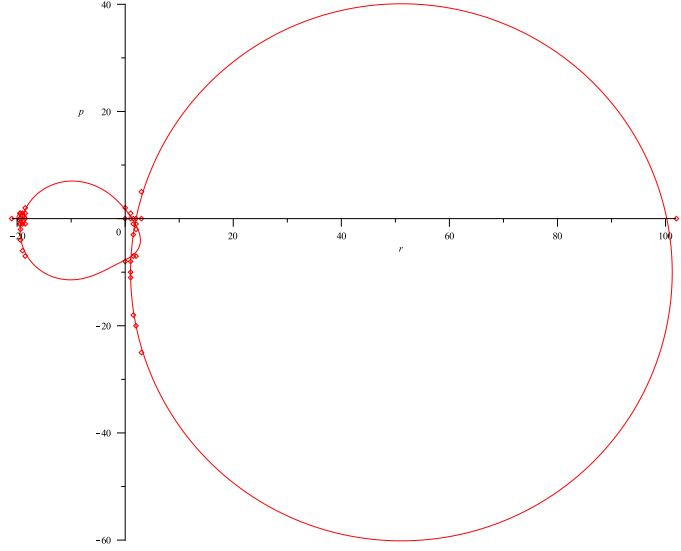}
\caption{Minimal algebraic curve containing border of $D$-decomposition and point cloud for an example 3.}
\end{figure}
But the border of $D$-decomposition is not an algebraic curve, because it does not contain part of the medium oval lying in the smallest one.
\begin{figure}\label{14.1}
\includegraphics[scale=0.3]{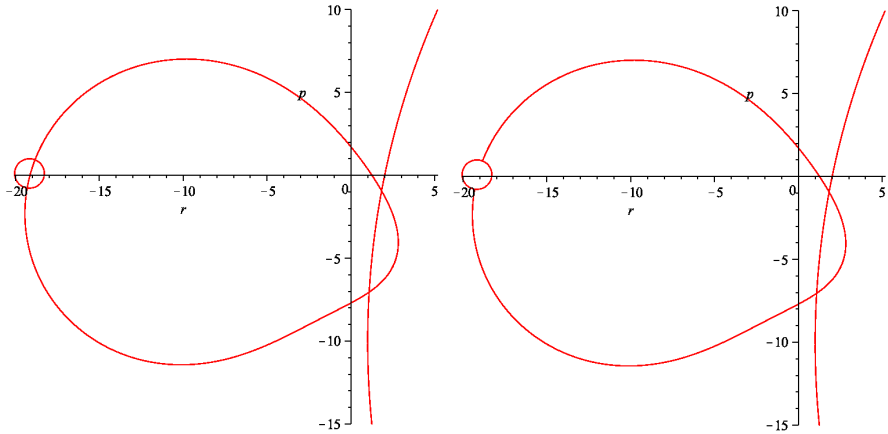}
\caption{Fragment of a minimal algebraic curve containing the border of $D$-decomposition for example 3(right side) and fragment of $D$-decomposition itself(left side).}
\end{figure}
\end{example}
\begin{example}
The next example is an output feedback of a form $K=\begin{pmatrix}-r&p\\p&r\end{pmatrix}$ for a co ntinuous-time systems defined by matrices
$$
\begin{array}{l}
A=\begin{pmatrix}79&20&-30&-20\\-41&-12&17&13\\167&40&-60&-38\\33,5&9&-14,5&-11\end{pmatrix},\\\\ B=\begin{pmatrix}0,219&0,9346\\0,047&0,3835\\0,6789&0,5194\\0,6793&0,831\end{pmatrix}, C=\begin{pmatrix}0,0346&0,5297&0,0077&0,0668\\0,0535&0,6711&0,3834&0,4175\end{pmatrix}.
\end{array}
$$
\begin{figure}\label{15}
\includegraphics[scale=0.6]{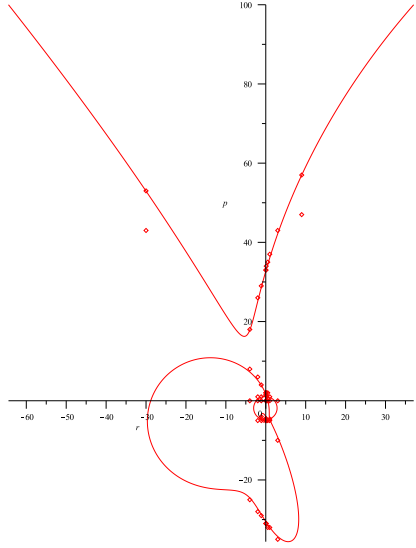}
\caption{Minimal algebraic curve containing the border of $D$-decomposition and point cloud representing all root invariant regions. Newly discovered upper component separates regions with 3 and 1 stable roots.}
\end{figure}

\begin{figure}\label{15.1}
\includegraphics[scale=0.3]{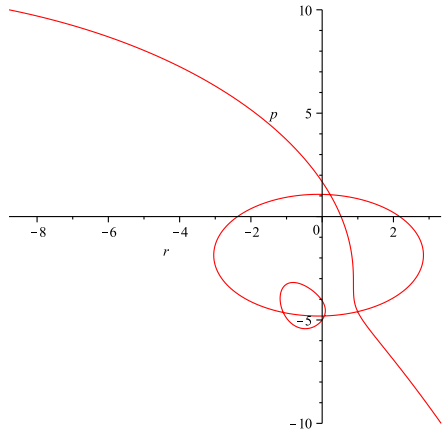}
\caption{Fragment of a minimal algebraic curve containing the border of $D$-decomposition for example 4.}
\end{figure}
This system has been considered in \cite{QBRDYD1995}, Example 2 and  \cite{GP2006}, Example 15.
We have found a new connected region of  $D$-decomposition with one stable root. This region isaway from all other components, its border corresponds to $\omega >20$.

This fact illustrates globality of our method. We don't need to prescribe any parameter ranges, but we do automatically get an answer for a whole space.
\end{example}
\begin{example}
Here we give an example of polynomial family with complex coefficient for which bound from Lemma \ref{lem:1} is sharp.

Namely, let
$$
\begin{array}{l}
P(s,r,p)=(4-22j)s^2+((30-20j)r-(145+8j)p+68)s+\\+((21-8j)r^2-(34+13j)p^2+(40+54j)rp+(2-26j)r+(-34+19j)p-(8+6j))
\end{array}
$$
Then the border of $D$-decomposition could be parametrised as:
$$
\begin{array}{l}
 4s^2+(20r+8p-4)s-21r^2+40rp+2r-21p^2-34p-8=0\\
-22s^2-145ps+68s+30rs-54rp-8r^2+19p-13p^2-26r-6=0
\end{array}
$$
An explicit equation of border of $D$-decomposition could be written as:
$$
\begin{array}{l}
1008000p^4-3303960p^3r+1782100p^2r^2+978760pr^3-627300r^4+62160p^3+\\+1366648p^2r+1219928pr^2-1200320r^3-1679328p^2+630352pr-145904r^2-\\-309120p+310272r+64512
\end{array}
$$
The decomposition  of plane into $11$ regions is shown on figure.
\begin{figure}\label{exp2}
\includegraphics[scale=0.3]{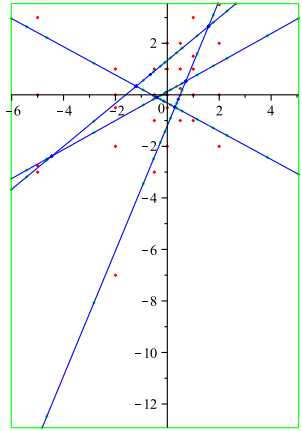}
\caption{The border of $D$-decomposition and point cloud for example 5.}
\end{figure}
\end{example}
\section{Conclusions}
We obtain a new technique of studying $D$-decomposition using real algebraic geometry. Method of obtaining an explicit equation for a border of $D$-decomposition is given as well as method of constructing point clouds from all root invariant regions.

An estimates for a number of root invariant regions are given in $n$-dimensional and in $2$-dimensional cases. The last bound is close to the sharp one.

This technique could be applied in construction of low order regulators as well as in robust stability analysis.

\end{document}